\documentclass[reqno, a4paper,12pt]{amsart} 
\usepackage{latexsym}
\usepackage[english]{babel}
\usepackage{fancyhdr}
\usepackage[mathscr]{eucal}
\usepackage{amsmath}
\usepackage{mathrsfs}
\usepackage{amsthm}
\usepackage{amsfonts}
\usepackage{amssymb}
\usepackage{amscd}
\usepackage{bbm}
\usepackage{graphicx}
\usepackage{graphics}
\usepackage{latexsym}
\usepackage{color}
\usepackage{ascmac}

\theoremstyle{plain}
\newtheorem{theorem}{Theorem}[section]
\newtheorem{remark}[theorem]{Remark}
\newtheorem{lemma}[theorem]{Lemma}
\newtheorem{corollary}[theorem]{Corollary}
\newtheorem{proposition}[theorem]{Proposition}
\newtheorem{assumption}[theorem]{Assumptions}

\theoremstyle{definition}
\newtheorem{definition}[theorem]{Definition}

\newcommand {\loc} {\textrm{loc}}

\def\Kg {{\mathcal K}}
\def\Xs {{\mathscr X}}
\def\Ys {{\mathscr Y}}
\def\Ts {{\mathscr T}}
\def\Ks {{\mathscr K}}
\def\tKg {{\tilde{\mathcal K}}}

\def\Wg {{\mathcal W}}
\def\Xg {{\mathcal X}}

\def\tM {{\tilde{M}}}
\numberwithin{equation}{section}

\def\ben{\begin{enumerate}}
\def\een{\end{enumerate}}
\def\bgdf{\begin{definition}}
\def\eddf{\end{definition}}
\def\bgpf{\begin{proof}}
\def\edpf{\end{proof}}
\def\bgth{\begin{theorem}}
\def\edth{\end{theorem}}
\def\bgprop{\begin{proposition}}
\def\edprop{\end{proposition}}
\def\bgrm{\begin{remark}}
\def\edrm{\end{remark}}
\def\bgcor{\begin{corollary}}
\def\edcor{\end{corollary}}
\def\bgass{\begin{assumption}}
\def\edass{\end{assumption}}

\def\lbeq(#1){\label{eqn:#1}}
\def\refeq(#1){{\textrm (\ref{eqn:#1})}}
\def\refeqs(#1,#2){{\textrm (\ref{eqn:#1}) and (\ref{eqn:#2})}}
\def\refeqss(#1,#2,#3)
{{\textrm (\ref{eqn:#1}),\ (\ref{eqn:#2}) and (\ref{eqn:#3})}}
\def\lbth(#1){\label{th:#1}}
\def\refth(#1){{\textrm Theorem \ref{th:#1}}}
\def\refths(#1,#2){{\textrm Theorems \ref{th:#1} and \ref{th:#2}}}
\def\refthss(#1,#2,#3){{\textrm Theorems \ref{th:#1}, \ref{th:#2} and 
\ref{th:#3}}}
\def\refthb(#1){{\textbf Theorem \ref{th:#1}}}
\def\lblm(#1){\label{lm:#1}}
\def\reflm(#1){{\textrm Lemma \ref{lm:#1}}}
\def\reflms(#1,#2){{\textrm Lemmas \ref{lm:#1} and \ref{lm:#2}}}
\def\reflmss(#1,#2,#3){{\textrm Lemmas \ref{lm:#1}, \ref{lm:#2} and \ref{lm:#3}}}
\def\reflmsss(#1,#2,#3,#4){{\textrm Lemmas \ref{lm:#1},\, \ref{lm:#2},\, \ref{lm:#3} and \ref{lm:#4}}}
\def\reflmb(#1){{\textbf Lemma \ref{lm:#1}}}
\def\lbprop(#1){\label{prp:#1}}
\def\refprop(#1){{\textrm Proposition \ref{prp:#1}}}
\def\refprops(#1,#2,#3,#4){{\textrm Propositions \ref{prp:#1},\, \ref{prp:#2},
\, \ref{prp:#3} \, and \ref{prp:#4}}}
\def\refpropb(#1){{\textbf Proposition \ref{prp:#1}.}}
\def\lbcor(#1){\label{cor:#1}}
\def\refcor(#1){{\textrm Corollary \ref{cor:#1}}}
\def\refcorb(#1){{\textbf Corollary \ref{cor:#1}}}
\def\refcors(#1,#2){{\textrm Corollaries \ref{cor:#1} and \ref{cor:#2}}}
\def\lbrm(#1){\label{rm:#1}}
\def\refrm(#1){{\textrm Remark \ref{rm:#1}}}
\def\refrmss(#1,#2,#3){{\textrm Remark \ref{rm:#1},\, \ref{rm:#2}\, and \ref{rm:#3}}}
\def\lbass(#1){\label{ass:#1}}
\def\refass(#1){{\textrm Assumption \ref{ass:#1}}}
\def\lbdf(#1){\label{df:#1}}
\def\refdf(#1){{\textrm Definition \ref{df:#1}}}
\def\refdfs(#1,#2){{\textrm Definitions \ref{def:#1} and \ref{def:#2}}}
\def\lbsec(#1){\label{s:#1}}
\def\refsec(#1){{\textrm \S\ref{s:#1}}}
\def\lbsubsec(#1){\label{ss:#1}}
\def\refsubsec(#1){{\textrm \S\ref{ss:#1}}}

\def\Ag{{\mathcal A}}
\def\Bg{{\mathcal B}}

\def\Kg{{\mathcal K}}

\newcommand{\lam}{\lambda}

\def\Bb{{\textbf B}}

\newcommand {\ph}{{\varphi}}

\def\bqn{\begin{equation}}
\def\eqn{\end{equation}}
\def\C{{\mathbb C}}
\def\N{{\mathbb N}}

\def\R{{\mathbb R}}
\def\Rg {{\mathcal R}}
\def\a{\alpha}

\def\Sg{{\mathcal S}}
\def\Ug{{\mathcal U}}

\def\p{\psi}

\def\ep{\varepsilon}
\def\z{\zeta}

\def\th{\theta}

\def\ka{\kappa}
\def\m{\mu}
\def\n{\nu}
\def\r{\rho}
\def\s{\sigma}
\def\t{\tau}

\def\W{\Omega}
\def\Hg {{\mathcal H}}
\def\Hs {{\mathscr H}}
\def\Ls {{\mathscr L}}
\def\Zg {{\mathcal Z}}
\def\la{\langle}
\def\ra{\rangle}

\def\lap{\Delta}

\def\pa{{\partial}}

\def\bglm{\begin{lemma}} 
\def\edlm{\end{lemma}} 
\def\br{\begin{array}}
\def\er{\end{array}}

\begin{document}

\title[N-body system with time dependent interactions]
{Boundedness of energy for N-body Schr\"odinger equations with time 
dependent small potentials}
\footnote{Supported by JSPS grant in aid for scientific research No. 19K03589}

\author[K.~Yajima]{Kenji Yajima}
\address{Department of Mathematics \\ Gakushuin University 
\\ 1-5-1 Mejiro \\ Toshima-ku \\ Tokyo 171-8588 (Japan)}. \\ 
\email{kenji.yajima@gakushuin.ac.jp}

\begin{abstract} 
We prove that Sobolev norms of solutions to 
time dependent Schr\"odinger equations for $d$-dimensional 
$N$-partcles interacting via time dependent two body 
potentials are bounded in time if certain Lebesgue norms of 
the potentials are small uniformly in time. The proof uses 
the scattering theory in the extended phase space 
which proves that all particles scatter freely 
in the remote past and far future. 

\end{abstract}

\date{}

\maketitle

\section{Introduction, Results} \label{sec:theorems} 

We consider Schr\"odinger equations for $N$ particles 
in $\R^d$, $d \geq 3$, 
\bqn
\lbeq(SE-1)
i \pa_t u(t,x)= \Big( -\sum_{j=1}^N \frac1{2m_j}\lap_{j} + 
\ka V(t,x)\Big) u(t,x)= \colon H_\kappa (t) u(t,x), 
\eqn
interacting via (complex) time dependent short range two-body potentials: 
\[   
V(t,x) = \sum_{1\leq j<k\leq N} V_{jk}(t,x_j-x_k) 
\] 
with {\it{small}} coupling constant $\ka$, where 
$x_j\in \R^{d}$ and $m_j$, $1\leq j \leq N$, are the 
position and the mass of $j$-th particle respectively, 
$x= (x_1, \dots, x_N)\in \R^{Nd}$, $\lap_{j}$ is the 
$d$-dimensional Laplacian with respect to $x_j$ and 
we have set $\hslash=1$. The purpose of the present paper 
is to show that the Sobolev norm $\|u(t,x)\|_{H^m(\R^{Nd})}$, 
$m=0,1,\dots$ of solutions \refeq(SE-1), hence the energy 
$(H_\ka(t)u(t), u(t))_{L^2(\R^{Nd})}$ of the system, 
is bounded in time $t \in \R$ if 
\bqn \lbeq(potential)
\sum_{|\a|\leq m}\sup_{t \in \R} 
\|\pa_x^\a V_{jk}(t,\cdot)\|_{(L^p\cap L^q)(\R^d)}\leq C  
\quad 1\leq j<k\leq N\,,
\eqn 
for some $p$ and $q$ such that $1\leq p<d/2<q \leq \infty$. 

If $V_{jk}(t,y)=V_{jk}(x)$, $1\leq j<k\leq N$ 
are real and independent of $t\in \R$,  
then $H_\ka$ on the right of \refeq(SE-1) with $\ka\in \R$ is selfadjoint in 
$\Hs= L^2(\R^{Nd})$ with domain $D(H_\ka)=H^2(\R^{Nd})$ and the 
energy $\la u(t), H_\ka u(t)\ra$ is conserved, hence, a fortiori 
bounded as $t \to \pm \infty$. However, if $V(t,x)$ is genuinly 
$t$-dependent, it is a suble question whether or not the energy of 
the system remains bounded in time. When $N=2$ and $V(t,x)$ 
is real, smooth and rapidly decreasing as 
$|x|\to \infty$ uniformly with respect to $t\in \R$,  
Bourgain (\cite{Bourgain}) has shown that Sobolev norms 
$\|u(t)\|_{H^m(\R^d)}$, $m \geq 0$, of solutions 
of \refeq(SE-1) remain bounded as $t\to \pm\infty$ 
if $\ka$ is sufficiently small and that, 
without the smallness condition, they satisfy  
$\|u(t)\|_{H^m(\R^d)}\leq C _\ep\la t\ra^\ep\|u(0)\|_{H^m(\R^d)}$ 
for any $\ep>0$, however, {\it the factor $\la t\ra^\ep$ cannot be removed 
when $m \geq 1$ in general}. Thus, we extend in this paper 
{\it the first part} of \cite{Bourgain} to 
$N$-body Hamiltonians with time dependent complex singular potentials. 

There are many works on the large $t$ behavior of 
Sobolev norms of the solutions of Schr\"odinger 
equations. However, except \cite{Bourgain} mentioned above, 
they deal with the case that the operator $H_{\ka}(t)$ on the right 
of \refeq(SE-1) has for 
each $t$ discrete spectrum  or $V(t,x)$ is periodic in $x$ 
and prove that Sobolev norms can increase only 
as slowly as $C |t|^\ep$ for arbitray small $\ep>0$ 
or $C\log |t|$ as $t \to \infty$ (see \cite{Bambusi}, 
\cite{Bour-2}, \cite{Bour-3}, \cite{Delort}, \cite{Duclos}, \
\cite{Erdogan}, \cite{Nenciu}, \cite{Wang} and references therein) 
and, as far as the author is aware of, there are no results so far 
for $N$-body system with genuinely time dependent potentials. 

For the reason to be explained below \refth(1), 
we consider \refeq(SE-1) 
for vector valued functions $u(t,x)\in \C^n$, $n=1,2, \dots$ 
with matrix potentials 
$V_{jk}(t,y)\in M(n)$:
\bqn \lbeq(SE-1a)
i\pa_t u = -\lap_X u  + {\ka}V(t,x) u, \quad 
V(t,x) = \sum_{1\leq j<k\leq N} V_{jk}(t,x_j-x_k)  
\eqn 
where $-\lap_X= \sum_{j=1}^N (2m_j)^{-1}\lap_j$ 
(see \refsec(existence)) and  $M(n) = M(n,n)$, 
$M(m,n)$ being the space of $m \times n$-matrices. 
For $M(m,n)$-valued function $A=(a_{jk})$ on a measure space $\W$, 
$A\in L^p(\W)$ means that $a_{jk}\in L^p(\W)$ for all $j,k$, 
$\|A\|_{L^p}=\left(\sum \|a_{jk}\|_{L^p}^p\right)^{\frac1{p}}$ 
and $\|A\|_p =\|A\|_{L^p}$. 
For Banach spaces $\Xs$ and $\Ys$ 
which are subspaces of a linear topological space $\Ts$, 
$\Xs \cap \Ys$ and $\Xs+ \Ys$ are Banach spaces with respective norms 
\[
\|u\|_{\Xs \cap \Ys}= \|u\|_{\Xs}+ \|u\|_{\Ys},\quad 
\|u\|_{\Xs+\Ys}
= \inf \{\|a\|_{\Xs}+ \|b\|_{\Ys}\colon u = a+ b \}; 
\]
$\Ls^p= L^p(\R^d:\C^n)$ 
and $\Ls_{p,q}= (L^{p}\cap L^q)(\R^d \colon M(n))$; 
$\dot A_{jk}(t,y)= \pa_t A_{jk}(t,y)$, etc. are $t$-derivatives. 

We remark before stating the theorem that, 
since $\Ls_{p,q} \subset L^d(\R^d \colon M(n))$ if $p<d<q$, 
$V_{jk}(t,y)$ of the following \refth(1) satisfy 
the conditions of \refthss(E,E-a,E-b) below with $p=d/2$ 
and \refeq(SE-1a) generates a unique  propagator $U(t,s, \kappa)$ on 
$\Hs \colon =L^2(\R^{Nd} \colon \C^n)$ (see \refdf(propagator)).  
\refth(1) holds for more general multi-particle 
interactions, however, we restrict ourselves 
to two-body interactions for notational simplicity.

\bgth \lbth(1) 
Let $n=1,2, \dots$ and $m\in \N \cup \{0\}$. 
Suppose that $V_{jk}(t,y)\in M(n)$, ${1\leq j<k\leq N}$, are factorized 
by $A_{jk}(t,y)\in M(n)$ and $B_{jk}(t,y)\in M(n)$,   
\[
V_{jk}(t,y)=B_{jk}(t,y)^\ast A_{jk}(t,y),
\]
which satisfy the following conditions for any $|\a|\leq m$:   
\ben
\item[\textrm (1)] There exist $p$ and $q$ such that 
$2\leq p < d <q \leq \infty$ and $q\geq 4$ if $d=3$ and such that
$\pa^\a_y A_{jk}(t,\cdot)$ and 
$\pa^\a_y B_{jk}(t,\cdot)$ are $\Ls_{p,q}$-valued functions of $t \in \R$ and  
\bqn 
\sup_{t \in \R} \left(
\|\pa^\a_y A_{jk}(t,\cdot)\|_{\Ls_{p,q}}+
\|\pa^\a_y B_{jk}(t,\cdot)\|_{\Ls_{p,q}}\right)<\infty\, .
\eqn 
\item[\textrm (2)] For a.e. $y \in \R^d$ they are 
locally absolutely continuous with respect to $t$ and 
\[
\pa^\a_y \dot {A}_{jk}(t,y), \pa^\a_y \dot B_{jk}(t,y)\in 
L^{2}_{\loc}(\R, (L^{r}+L^\infty)(\R^d \colon M(n)))
\] 
for $r=2$ if $d =3$ and $r=d/2$ if $d \geq 4$. 
\een 
Then, for $|\ka|<\ka_m$, $\ka_m>0$ being a small constant, 
there exits a $C_{\ka,m}<\infty$ such that 
\bqn \lbeq(0-3)
\sup_{t,s\in \R}\|U(t,s,\ka)\ph\|_{H^m(\R^{Nd} \colon \C^n)} 
\leq C_{\ka,m} 
\|\ph\|_{H^m(\R^{Nd}\colon \C^n)}.
\eqn 
\edth 

For the proof of \refth(1) we use the scattering theory for 
\refeq(SE-1) and the induction argument 
of Bourgain(\cite{Bourgain}) that the boundedness of 
$\|u(t,\cdot)\|_{H^{k+1}(\R^d)}$ 
follows from that of $\|{\bf u }(t,\cdot)|_{H^k(\R^d\colon \C^{(Nd+1)})}$  
of ${\textbf u}= \begin{pmatrix} u \\ \nabla_x u \end{pmatrix}$ 
which satisfies  
\bqn \lbeq(syst-0)
i \pa_t {\textbf u}= H_0 \textbf{u}+ \ka 
\begin{pmatrix} V(t,x)  &  {\textbf 0} \\ 
{\nabla_x}V (t,x) & V(t,x){\textbf 1} \end{pmatrix}
{\textbf u}\,.
\eqn 
In \refsec(existence) we recall the results of \cite{Ya-2016} 
which imply that \refeq(SE-1a) generates a unique propagator 
$\{U(t,s,\ka)\}$ on $\Hs$. Then, in \refsec(howland), 
following Howland \cite{Howland},
we introduce the extended phase space 
$\Ks = L^2(\R) \otimes \Hs$ and the strongly continuous one parameter group
of bounded operators $\{\Ug(\s) \colon \s\in \R\}$ 
on $\Ks$ by $(\Ug(s)u)(t) = U(t, t-s) u(t-s)$. 
Let $(\Ug_0(\s)u)(t)= e^{-i{\s}H_0}u(t-s)$. We then state \refth(2) 
in \refsec(wave) that strong limit   
$\Wg_{+}= \lim_{\s\to \infty}e^{i\s\Kg}e^{-i\s\Kg_0}$ exists in $\Ks$,  
it is an isomorphism of $\Ks$ and satisfies 
the intertwing property $e^{-i\s\Kg}= \Wg_{+} e^{-i\s\Kg_0}\Wg_{+}^{-1}$. 
Postponing the proof of \refth(2) to \refsec(proof-2), we show 
also in \refsec(wave) that \refth(2) implies that the limit in $\Hg$ 
\[
\lim_{t\to \infty} U(s,t,\ka)e^{-i(t-s)H_0} = W_{+}(s)
\]
exists, it is an isomorphism of $\Hs$ and is  
uniformly bounded for $s\in \R$ along with the inverse $W_{+}(s)^{-1}$:
\bqn \lbeq(unif-bdd)
\|W_{+}(s)\|_{\Bb(\Hs)}\leq C, \quad 
\|W_{+}(s)^{-1}\|_{\Bb(\Hs)}\leq C, \quad s \in \R
\eqn 
and that it satisfies the intertwing property:
\bqn \lbeq(interwine)
U(t,s) = W_{+}(t)e^{-i(t-s)H_0} W_{+}(s)^{-1}, \quad t,s \in \R.
\eqn 
We prove \refth(1) in \refsec(proof-1). The case $m=0$ 
is evident from \refeq(unif-bdd) and \refeq(interwine). 
If we assume that \refth(1) holds for $m=1, \dots, k$, then 
the ''potential'' of \refeq(syst-0) satisfies the condition of \refth(1) 
for $m=k$ as $M((Nd+1)n)$-valued function. 
Hence, $\|{\textbf u}(t)\|_{H^k(\R^{Nd}: \C^{(Nd+1)n)})}$ 
is bounded for $t \in \R$ by the induction hypothesis.
and $\|u(t)\|_{H^{k+1}(R^{Nd}:\C^n)}\leq C$. 
The proof of \refth(3) is given in \refsec(proof-2) by adapting 
Ioirio-O'Carrol's argument (\cite{IO}) to 
the extended phase space and by applying Kato's theory of smooth 
perturbations (\cite{Kato-e}). Thus, the method employed here 
are rather old theory of scattering and it is expected that  
its modern theory can produce more refined result.

\section{Existence and regularity of propagators} \lbsec(existence)
We begin with recalling the result of \cite{Ya-2016} on the existence and 
the regularity of the propagator for \refeq(SE-1a) in the form modified 
for our purpose. 
We use some $N$-body notation due to Agmon \cite{Agmon}: $X$ is the space  
$\R^{Nd}$ with the so called {\it mass inner product}
\[ 
(x,y)_m \colon = \sum_{j=1}^N 2m_j (x_j, y_j)_{\R^d}, 
\quad x= (x_1, \dots, x_N), \ y= (y_1, \dots, y_N), 
\]
then $\sum_{j=1}^N (2m_j)^{-1}\lap_{j}= \lap_X$;  
for the pair $\{j,k\}$, $1\leq j,k\leq N$, we set   
\begin{gather*}
X_{jk}^{out} =\{x \in X \colon x_j= x_k=0\}, \quad 
X_{jk}^{in} = X \ominus X_{jk}^{out} \simeq \R^{2d}, \\
X_{jk}^c= \{x\in X_{jk}^{in} \colon x_j=x_k\}\simeq \R^{d}, \quad 
X_{jk}^r= X_{jk}^{in}\ominus X_{jk}^c \simeq \R^{d}. 
\end{gather*}
Then, 
$X= X_{jk}^{out} \oplus X^{in}_{jk}$ and 
$X^{in}_{jk}= X_{jk}^c\oplus X_{jk}^r$.  We let 
\bqn 
X_{jk}= X_{jk}^{out} \oplus X_{jk}^c
\eqn 
and the corresponding orthogonal decomposition of $x\in X$ be 
\[
x= x^{jk} \oplus x_{jk}^c \oplus x_{jk}^{out} \in X_{jk}^r \oplus X^c_{jk} 
\oplus X_{jk}^{out}; \quad x_{jk}=\colon x_{jk}^c \oplus x_{jk}^{out}\in X_{jk}. 
\]
This orthogonal decomposition leads to 
\begin{gather} 
\Hs = L^2(X_{jk}^r \colon \Hs_{jk}), \quad  \Hs_{jk}\colon = 
L^2(X_{jk}\colon \C^n) .  \\
\lap_x = \lap_{x^{jk}}\otimes \textbf{1}_{\Hs_{jk}}+ 
\textbf{1}_{L^2(X_{jk}^r)} \otimes \lap_{x_{jk}}.   \lbeq(lap-decompo)
\end{gather} 
We denote $\lap_{jk}= \lap_{x^{jk}}$. It is easy to check that  
\begin{align*}
x^{12}& = \left(x^{12}_1, x^{12}_2 , 0, \dots, 0\right), 
\ \ 
x^{12}_1/m_2= -x^{12}_2/ m_1= (x_1-x_2)/(m_1+m_2), 
\\
x_{12}^c& = \left(x_{12}^c,x_{21}^c, 
0, \dots, 0 \right), \quad x_{12}^c=(m_1x_1+m_2x_2)/(m_1+m_2)
\end{align*}
and $x_{12}^{out}= (0,0,x_3, \dots, x_N)$; 
we have similar identities for $x^{jk}$, $x_{jk}^c$ and $x_{jk}^{out}$ 
and $V_{jk}(t,x_j-x_k)$ is a function of $(t,x^{jk})$. Abusing notation, 
we often write $V_{jk}(t,x^{jk})$ for $V_{jk}(t,x_j-x_k)$.

\bgdf \lbdf(propagator) 
We say a strongly continuous family 
$\{U(t,s,\ka) \colon s, t\in \R\}$ of bounded operators 
on $\Hs= L^2(X:\C^n)$ 
is the propagator for \refeq(SE-1a) if it satisfies 
the Champan-Kolmogorov equation: $U(t,s,\ka)U(s,r,\ka)=U(t,r,\ka)$, 
$U(t,t,\ka)= {\textbf 1}_{\Hs}$ for $t,s,r\in \R$ and, 
if $u(t,x)=(U(t,s,\ka)\ph)(x)$, $\ph \in \Hs$, 
is a solution of \refeq(SE-1) such that $u(s,x)= \ph(x)$. 
\eddf 

For a given 
$d/2 \leq p \leq \infty$ and compact intervals $I$, we let 
\bqn \lbeq(X)
\Xg^p(I)= \bigcap\limits_{1\leq j<k\leq \infty}L^{\th}(I, 
L^{\ell}(X_{jk}^r, \Hs_{jk}))\,;
\eqn 
where $\th= {4p}/{d}(\geq 2)$ and ${\ell}=2p/(p-1)(\geq 2d/(d-2))$; 
$\Xg^p_{\loc}$ denotes the space of $f$ such that $f \in \Xg^p(I)$ 
for compact intervals. We remark that the uniqueness 
of the solution in the following theorem is under the condition 
$u(t, \cdot) \in C(\R, \Hs) \cap \Xg^p_{\loc}$. 

\bgth \lbth(E) 
Let $\{V_{jk}(t,y)\}_{1\leq j<k\leq N}$ be $M(n)$-valued functions 
such that 
$V_{jk}(t,y)\in C(\R, L^{p}(\R^d))+ C(\R, L^\infty(\R^d))$ 
for a $d/2 \leq p \leq \infty$ and 
\bqn 
\sup_{t\in \R} (\|V_{jk}(t)\|_{L^p(\R^d)} + \|V_{jk}(t)\|_{L^\infty(\R^d)})
\leq C <\infty. 
\lbeq(C-1)
\eqn
Then, Eqn. \refeq(SE-1a) uniquely 
generates a propagator $U(t,s,\ka)$ such that for evey  $\ph \in \Hs$ 
$u(t,x)\colon =(U(t,s,\ka)\ph)(x)\in C(\R, \Hs) \cap \Xg^p_{\loc}$. 
We have $u(t, \cdot) \in C^1(\R, H^{-2}(X:\C^n))$ and \refeq(SE-1a) 
is satisfied in $H^{-2}(X, \C^n)$. Moreover:
\ben
\item[\textrm{(1)}] For any $a>0$, 
$\sup_{t,s\in \R\colon |t-s|\leq a}\|U(t, s,\ka)\|_{\Bb(\Hs)}=C_a<\infty$ 
and 
\bqn \lbeq(local-Strichartz)
\sup_{s\in \R} \|u(t)\|_{\Xg^p([s-a,s+a])} \leq C_a \|\ph\|_{\Hs}\,.
\eqn 
\item[\textrm{(2)}] $U(t,s,\ka)$  is unitary if all 
$V_{jk}(t,y)$ are Hermitian and $\ka\in \R$. 
\een
\edth 

The second theorem is on the regularity of the solution obtained 
in \refth(E).   

\bgth \lbth(E-a) Suppose $\{V_{jk}(t,y)\}_{1\leq j<k\leq N}$ 
satisfy \refeq(C-1) with $\tilde{p}=\max(2,p)$ replacing $p$ 
and, in addition, 
\begin{gather}
\dot V_{jk}(t,y) \in 
L^{b}_{\loc}(\R, L^{q}(\R^d)) + L^1_{\textrm{loc}}(\R, L^\infty(\R^d)),  
\lbeq(C-2) \\
1/b= 1-d/4p; \ 1/q=\left\{ \br{l}1/2+1/2p \ \mbox{if} \ d =3, \\ 
2/d+ 1/2p \ \mbox{if}\ d\geq 4. \er \right.  \lbeq(C-2a)
\end{gather} 
Then, for $\ph \in H^2(X\colon \C^n)$,  
$u(t,x)=U(t,s,\kappa)\ph$ of \refth(E) satisfies 
\bqn \lbeq(E-a)
u(t,x) \in C^1(\R, \Hs)
\cap C(\R, H^2(X: \C^n)), \quad  
\dot u (t,x)\in \Xg^{{p}}_{\loc}.
\eqn 
\edth 

We remark that $\tilde{p}=p$ when $d\geq 4$; if $d=3$ and $p<2$ and if 
$V_{jk}(t,y)$ satisfies \refeq(C-1) with $\tilde{p}=\max(2,p)$ 
replacing $p$, then it also satisfies \refeq(C-1) with the orginal $p$.

\refths(E,E-a) except statement (1) of \refth(E) are stated and 
proved in \cite{Ya-2016} for the case $n=1$ with real valued 
$V_{jk}(t,y)$ and, with $\Sigma(2)$  
and $\Sigma(-2)$ 
in place of Sobolev spaces $H^2(X)$ and $H^{-2}(X)$ respectively, 
where $\Sigma(2)=\{u \colon \|u\|_{\Sigma(2)}^2= 
\sum_{|\alpha|+|\beta|\leq 2}\|x^\alpha \pa^\beta u\|_{\Hs}^2<\infty\}$   
is domain of the harmonic oscillator $-\lap + x^2$ and $\Sigma(-2)$ 
is its dual space. However, the extension to $n \geq 2$ 
with matrix-valued $V_{jk}(t,y)$ is obvious; 
if external electro-magnetic fields are absent, then 
the propagator $e^{-it\lap}$ for the independent particles has 
Sobolev spaces $H^m(X: \C^n)$ as invariant subspaces and, 
$\Sigma(2)$ and $\Sigma(-2)$ 
may be replaced by $H^2(X)$ and $H^{-2}(X)$ respectively. 
Statement (1) of \refth(E) is also evident since 
assumption \refeq(C-1) 
is translation invariant with respect to $t \in \R$.
We omit the proof of \refths(E,E-a),leaving for the readers 
to check the details. 

If $\ph \in H^2(X,\C^n)$ and $u(t,x)= (U(t,s,\ka)\ph)(x)$, then it is evident 
that 
$\textbf{u}(t,x)=\begin{pmatrix} u(t,x) \\ \nabla_x u(t,x) 
\end{pmatrix}\in C^{(1+Nd)n}$ satisfies \refeq(syst-0).  
The next theorem implies that $\textbf{u}(t,x)$ is 
a unique solution of \refeq(syst-0).  For shortening formulas we write 
$\Hs^1=H^1(X\colon \C^{n(1+Nd)})$. 

\bgth \lbth(E-b) Let $\{V_{jk}(t,y)\}$ and 
$\{\nabla_y V_{jk}(t,y)\}$ satisfy \refeq(C-1) with $p$ 
being replaced by $\tilde{p}=\max(2,p)$ and  $\ph \in H^2(X, \C^n)$. 
Then, $\textbf{u}(t,x)$ satisfies \refeq(syst-0). If $\{V_{jk}(t,y)\}$ 
in addition satisfy \refeq(C-2),  
then ${\textbf u}$ is a unique solution of \refeq(syst-0) with 
the initial condition  
${\textbf u}(s,x)= \begin{pmatrix} \ph(x) \\ \nabla_x \ph(x)\end{pmatrix}$ 
and ${\textbf u} \in C(\R, \Hs^1) (\subset C(\R, \Hg) \cap \Xg^p_{\loc})$.  
\edth 
\bgpf If $\ph \in H^2(X, \C^n)$, then \refth(E-a) and 
the Sobolev embedding theorem imply that 
${\textbf u}\in C(\R, H^1(\R^d)) \subset C(I, \Hg) \cap \Xg^p_{\loc}$.  
Then, by \refth(E), ${\textbf u}(t,x)$ is the unique solution of 
\refeq(syst-0) which satisfies the condition of the theorem. 
\edpf 

\section{Howland scheme} \lbsec(howland) 

For proving \refth(1) we use the scattering theory for time dependent 
potentials. Following Howland \cite{Howland}, we introduce 
the extended phase space by $\Ks = L^2(\R^1, \Hs)= L^2(\R^1) \otimes \Hs$, 
where $\Hs= L^2(X \colon \C^n)$ in what follows,  
and ``the free propagator'' on $\Ks$ by 
\[
(\Ug_0(\s)u)(t)= e^{-i\s H_0}u(t-\s), \quad t \in \R\,, 
\] 
where $H_0= -(\lap_X) \otimes {\textbf 1}_{\C^n}$. 
It is easy to see that $\Ug_0(\s)$ is the unitary group 
generated by the selfajoint operator $\Kg_0$ on $\Ks$:
\bqn \lbeq(U0)
\Kg_0 = (-i\pa_t) \otimes {\textbf 1}_{\Hs} + 
{\textbf 1}_{L^2(\R^1)}\otimes H_0.
\eqn 

When condition (1) of \refth(1) is satisfied for $m=0$, then 
$\{V_{jk}(t,y)\}$ satisfy 
\refeq(C-1) and eqn. \refeq(SE-1a) 
uniquely generates a propagator $\{U(t,s,\ka)\}$. 
We define a one-parameter family of bounded of operators 
$\{\Ug(\s,\ka) \colon \s \in \R\}$ on $\Ks$ by  
\bqn \lbeq(Howland-1)
(\Ug(\s,\ka)u)(t)= U(t,t-\s,\ka)u(t-\s,\ka), \quad t \in \R\,. 
\eqn
In what follows untill the last section, assuming 
$|\ka|<\ka_0$ for a sufficiently small $\ka_0$, 
we omit the coupling constant $\ka$ 
from various quantities 

\bglm[\cite{Howland}] 
The family $\{\Ug(\s) \colon \s\in \R\}$  is a strongly continuous 
one-parameter group of bounded operators on $\Ks$. 
For constants $C\geq 1$ and $M\geq 0$ it satisfies    
\bqn \lbeq(bdd)
\|\Ug(\s)\|_{\Bb(\Ks)} \leq C e^{|\s|M}, \quad \s \in \R. 
\eqn 
There exits a unique closed operator $\Kg$ such that 
\bqn \lbeq(generator)
\Ug(\s)= e^{-i\s\Kg}, \quad -\infty<\s<\infty, 
\eqn 
the spectrum of $\Kg$ lies in the strip $\{\z \colon |\Im \z| \leq M \}$
and 
\bqn \lbeq(Phillips)
(\Kg -\z)^{-1} = \pm i \int_0^\infty e^{\mp i(\Kg-\z)\s} d\s, 
\quad \pm \Im \z >M.
\eqn 
\edlm 

\bgrm It will be proved in the next subsection 
that $M=0$ and ``$\Kg= \Kg_0 + V$''.
\edrm 

\bgpf The Chapman-Kolmogorov identity implies 
$\Ug(\s)\Ug(\t)= \Ug(\s+\t)$ for $\s, \t \in \R$. By \refth(E) (1) 
we have for an arbitrary $\s_0$ that      
\[
\|\Ug(\s)u\|_{\Ks}^2 = \int_{\R}\|U(t, t-\s)u(t-\s)\|_{\Hs}^2 dt 
\leq C_{\s_0} \|u\|_{\Ks}^2, \quad |\s|<\s_0.
\]
It follows that 
$\{\Ug(\s) \colon \s\in \R\}$ defines a one-parameter 
group of bounded operators on $\Ks$ which is locally uniformly bounded. 
Let $u \in  C_0(\R,\Hs)$ be supported by the interval $[a,b]\Subset \R$.
Then, for $|\s|\leq \s_0$, 
$\|U(t,t-\s)u(t-\s)- u(t)\|_{\Hs}\in C_0([a-\s_0, b+\s_0])$ 
and it converges to $0$ everywhere for $t \in \R$ as $\s\to 0$. Hence  
\[
\|\Ug(\s)u- u\|_{\Ks}^2 = \int_{\R} \|U(t,t-\s)u(t-\s)- u(t)\|_{\Hs}^2 dt \to 0 
\quad (\s \to 0). 
\]
Since $C_0(\R,\Hs)$ is dense in $\Ks$, this implies $\Ug(\s)$ is 
strongly continuous and $\{\Ug(\s)\}$ is a 
$C_0$-group of bounded operators. Then, the rest 
is well known, see e.g. \cite{Yosida}.  
\edpf

\section{Wave operators}  \lbsec(wave) 

In \cite{Bourgain}, the boundedness 
of $\|u(t)\|_{H^s(X)}$ for the case $N=2$ is proved by using 
the local smoothing property and the 
dispersive estimate for \refeq(SE-1). Here we shall prove it 
for $N =2,3, \dots$ via the following \refth(2) 
on wave operators of scattering theory which will be proved 
in \refsec(proof-2). 

\bgth \lbth(2) Let $n=1,2, \dots$. 
Suppose that $\{V_{jk}(t,y)\}_{1\leq j<k\leq N}\subset M(n)$ 
are factorized as 
$V_{jk}(t,y)=B_{jk}^\ast (t,y)A_{jk}(t,y)$ by 
$A_{jk}(t,y), B_{jk}(t,y)\in M(n)$ which satisfy 
condition (1) of \refth(1) for $m=0$. 
Then there exists a $\ka_0$ such that, for all $|\ka|<\ka_0$, the 
following statments are satisfied:
\ben 
\item[{\textrm (1)}] Both of the following strong limits in $\Ks$ exist: 
\begin{align} 
& \Wg_\pm (\ka) = \lim_{\s \to \pm \infty} \Ug (-\s, \ka) \Ug_0(\s), \lbeq(W)\\
& \Zg_\pm (\ka) = \lim_{\s \to \pm \infty} \Ug_0 (-\s) \Ug(\s,\ka)\,.  \lbeq(Z)
\end{align} 
They satisfy the following identities and are isomorphisms in $\Ks$: 
\bqn \lbeq(isom-H)
\Wg_\pm (\ka)\Zg_\pm(\ka)= \Zg_\pm(\ka)\Wg_\pm (\ka)={\textbf 1}_{\Ks}\,.
\eqn
\item[{\textrm (2)}] The operators $\Wg_\pm (\ka)$ satisfy 
the intertwing property: 
\bqn \lbeq(inter)
\Ug (\s, \ka)= \Wg_\pm (\ka)\Ug_0(\s)\Wg_\pm (\ka)^{-1}.  
\eqn  
 \een
\edth 

We postpone the proof of \refth(2) and proceed to 
proving \refth(1), admitting \refth(2) has already been proved. 
We first prove that \refth(2) implies the existence and the 
completeness of the wave operators for \refeq(SE-1a).

\bgth \lbth(3) Let $n=1,2, \dots$. 
Suppose the condition of \refth(2) is satisfied 
and $|\ka|$ is sufficiently small. Then, the following statements are 
satisfied: 
\ben 
\item[{\textrm (1)}] Both of the strong limits  
\begin{align}
& W_\pm (s,\ka)= \lim_{t \to \pm \infty} U(s,t,\ka)e^{-i(t-s)H_0}, 
\lbeq(Wa) \\
& Z_\pm (s,\ka) = \lim_{t \to \pm \infty} e^{i(t-s)H_0}U(t,s,\ka) \lbeq(Za)
\end{align}
exist in $\Bb(\Hs)$ for all $s\in \R$. They satisfy 
\bqn \lbeq(isom)
W_\pm (s,\ka)Z_\pm (s,\ka)= 
Z_\pm (s,\ka)W_\pm (s,\ka)= {\textbf 1}_{\Hs} 
\eqn 
and are isomorphisms in $\Hs$. 
\item[{\textrm (2)}] There exists an $s$-independent constant $C$ such that 
\bqn 
\|W_\pm (s,\ka)\|_{\Bb(\Hs)}\leq C, \quad 
\|Z_\pm (s,\ka)\|_{\Bb(\Hs)} \leq C .
\eqn  
\item[{\textrm (3)}] 
The wave operators satisfy the intertwing property:
\bqn \lbeq(inter-a)
U(t,s,\ka) = W_\pm(t,\ka)e^{-i(t-s)H_0} Z_\pm(s,\ka)\,.
\eqn
\een
\edth 

\bgrm Statement {\textrm (1)} of \refth(3) implies that for any $\ph \in \Hs$, 
the solution $U(t,s,\ka) \ph$ becomes asymptotically free as $t \to \pm \infty$:
\[
\lim_{t \to \pm \infty}\|U(t,s,\ka) \ph - e^{-i(t-s)H_0} \ph_\pm\|_{\Hs}=0, \quad 
\ph_\pm = Z_\pm (s,\ka)\ph.
\]
\edrm 

\bgpf  We prove the $+$ case only. 
The proof for the other case is similar. 
\refth(2) implies that for any $u \in \Ks$ the limit 
\bqn \lbeq(limit-1)
\lim_{\s\to\infty} \Ug(-\s) \Ug_0(\s) u(t) 
=\lim_{\s\to\infty} U(t,t+\s) U_0(t+\s,t)u(t)= (\Wg_{+} u)(t)
\eqn 
exists in $\Ks$. Let $u_0 \in \Hs$, $a>0$ be arbitrary and set  
$u(t) = U_0(t,0)u_0$ for $-a \leq t \leq a$ and $u(t)=0$ for 
$|t|>a$. Then, $u \in \Ks$ and, on substituting this $u(t)$ in 
\refeq(limit-1), we obtain that  as $\s \to \infty$ 
\[
\int_{-a}^a \|U(t,\s) U_0(\s,0)u_0 - (\Wg_{+} u)(t)\|_\Hs^2 dt \to 0 
\quad (\s \to \infty).
\] 
Since $\|U(0,t)\|_{\Bb(\Hs)}\leq C$ for $t \in [-a,a]$, it follows that 
\[
\lim_{\s\to \infty}
\int_{-a}^a \|U(0,\s) U_0(\s,0)u_0 - U(0,t)(\Wg_{+}u)(t)\|_\Hs^2=0\,.
\] 
Then Schwarz's inequality implies that, as $\s \to \infty$ ,  
\begin{multline*}
\Big\|
\int_{-a}^a (U(0,\s) U_0(\s,0)u_0 - U(0,t)(\Wg_{+} u)(t)) dt \Big\|_{\Hs}\\
\leq \sqrt{2a}
\Big(
\int_{-a}^a \|U(0,\s) U_0(\s,0)u_0 
- U(0,t)(\Wg_{+} u)(t)\|_{\Hs}^2 dt\Big)^\frac12 
\to 0. 
\end{multline*}
Thus, $U(0,\s) U_0(\s,0)u_0$ converges in $\Hs$ as $\s \to \infty$. 
Since $u_0 \in \Hs$ is arbitrary, this implies that 
for any $t$, the strong limit of \refeq(Wa):
\[
W_{+}(t)= \lim_{\s \to \infty} U(t,\s) U_0(\s,t)= 
U(t,0) (\lim_{\s \to \infty} U(0,\s) U_0(\s,0)) U_0 (0,t)= 
\]
exists in $\Hs$ for any $t\in \R$ and \refeq(limit-1) implies 
\bqn \lbeq(W=Wg)
(\Wg_{+} u)(t) = W_{+}(t)u(t), \quad a.e. \ t\in \R.
\eqn 
The proof for $Z_{+}(s)$ is similar.  They are bounded operators 
in $\Hs$ as they are strong limits of bounded operators. We then have 
\[
W_{+}(s)Z_{+}(s)=\lim_{t \to \infty} 
U(s,t,\ka)e^{-i(t-s)H_0}e^{i(t-s)H_0}U(t,s,\ka)=\textbf{1}  
\]
and, likewise, $Z_{+}(s)W_{+}(s)=\textbf{1} $, which imply 
the identity \refeq(isom) 
and $W_{+}(s)$ and $Z_{+} (s)$ are isomorphisms in $\Hs$. 

\noindent 
(2) The Chappman-Kolmogorov identity and the strong convergence imply  
\[
W_{+}(s) u_0 =U(s,t)W_{+}(t) U_0(t,s)u_0 
\]
for any $t,s\in \R$ and any $u_0 \in \Hs$. Fix $s\in \R$ and $u_0 \in \Hs$ 
arbitary and define $u(t) = U_0(t,s)u_0$ 
for $|t-s|\leq 1$ and $u(t) = 0$ for $|t-s|>1$. 
Then, since $\|U(t,s)\|_{\Bb(\Hs)}\leq C$ for $|t-s|\leq 1$ 
and $\|U_0(t,s)\|_{\Bb(\Hs)}=1$, 
\[
\|W_{+}(s) u_0\|_{\Hs} 
\leq  \frac1{2}\int_{s-1}^{s+1}\|U(s,t)W_{+}(t)u(t)\|_{\Hs}dt 
\leq  \frac{C}{\sqrt{2}}
\left(\int_{s-1}^{s+1}\|W_{+}(t)u(t)\|_{\Hs}^2dt\right)^\frac12 
\]
and, by using \refeq(W=Wg), we estimate the right side by 
\[
\frac{C\|\Wg_{+}\|_{\Bb(\Ks)}}{\sqrt{2}}
\left(\int_{s-1}^{s+1}\|u(t)\|_{\Hs}^2dt\right)^\frac12 
\leq C \|\Wg_{+}\|_{\Bb(\Ks)}\|u_0\|\,.
\]
Here $C$ is independent of $s\in \R$ and hence, 
$\|W_{+}(s)\|_{\Bb(\Hs)}$ is uniformly bounded for $s\in \R$. 
Likewise we obtain $\|Z_{+}(s)\|_{\Bb(\Hs)}\leq C$. 

\noindent 
(3) The standard argument for the selfadjoint case applies for 
\refeq(inter-a). We omit the details. 
\edpf

\section{Proof of \refthb(1)} \lbsec(proof-1)

We prove \refth(1) by induction on $m$ admitting \refth(3). 

\paragraph{\textbf Proof of \refthb(1) for the case $m=0$ 
and $n=1,2, \dots$}. 
Since $\{W_\pm(t)\colon t\in \R\}$ and 
$\{Z_\pm(s) \colon s\in \R\}$ are uniformly bounded in $\Hs$ and 
$e^{-itH_0}$ is unitary for $t \in \R$, we have from \refeq(inter-a) 
that for $\ph \in \Hs$  
\[
\|U(t,s)\ph\|_{\Hs} = \|W_\pm(t)e^{-i(t-s)H_0} Z_\pm(s)\ph\|_{\Hs}
\leq C \|\ph\|_{\Hs}
\]
for a constant $C>0$ independent of $(t,s)$. This proves 
\refth(1) for the case $m=0$. 

\paragraph{\textbf Proof of \refthb(1) for general $m=1,2, \dots,$}
We assume that \refth(1) has already been proved for $m=0, \dots, k$ 
and for all $n=1,2, \dots$ and suppose $m=k+1$. 
We prove 
\bqn \lbeq(0-3aa)
\|U(t,0)\ph\|_{H^{k+1}(X\colon \C^n)} \leq C_{k+1} 
\|\ph\|_{H^{k+1}(X \colon \C^n)}, \ n=1,2, \dots.
\eqn
For shortening formulas we write $\Hs^s = H^s(X\colon \C^{n(Nd+1)})$ 
and $U(t)$ for $U(t,0)$. It suffices \refeq(0-3aa) for 
$\ph \in H^{k+2}(X\colon \C^n)$ which is dense in $H^{k+1}(X \colon \C^n)$.
Then, since $k+2 \geq2$, \refth(E-a) implies 
\[
{\textbf u}(t)\colon = \begin{pmatrix} U(t)\ph \\  \nabla_x U(t)\ph \end{pmatrix}
\in C(\R, \Hs^1) \cap C(\R, \Hs^0)\subset \Xs^p_{\loc}
\]
is the unique solution of the initial value 
problem for the Schr\"odinger equation  \refeq(syst-0):  
\bqn \lbeq(5-1)
i \pa_t {\textbf u} 
= -\lap_x {\textbf u}+ 
\ka \sum_{1\leq j<k \leq N} 
\begin{pmatrix} V_{jk}(t)  &  {\textbf 0} \\ 
{\nabla}V_{jk}(t) & V_{jk}(t){\textbf 1}\end{pmatrix} {\textbf u}, \  
{\textbf u}(0)\colon = \begin{pmatrix} \ph \\  \nabla_x \ph \end{pmatrix}
\eqn 
and, the ``potential'' in \refeq(5-1) may be factorized as 
\[
\begin{pmatrix} V_{jk}(t)  &  {\textbf 0} \\
{\nabla}V_{jk}(t) & V_{jk}(t){\textbf 1}\end{pmatrix}
= \begin{pmatrix} B_{jk}(t)^\ast  &  ({\nabla}B_{jk}(t))^\ast \\
{\textbf 0} & B_{jk}(t)^\ast {\textbf 1}\end{pmatrix}^{\ast} 
\begin{pmatrix} A_{jk}(t)  &  {\textbf 0} \\
{\nabla}A_{jk}(t) &  A_{jk}(t){\textbf 1}\end{pmatrix}
\]
and the factors 
\[
A_{jk}'(t,y)= \begin{pmatrix} A_{jk}(t,y)  &  {\textbf 0} \\
{\nabla}A_{jk}(t,y) &  A_{jk}(t,y){\textbf 1}\end{pmatrix}, 
\ 
B_{jk}'(t,y)= \begin{pmatrix} B_{jk}(t,y)^\ast  &  ({\nabla}B_{jk}(t,y))^\ast \\
{\textbf 0} & B_{jk}(t,y)^\ast {\textbf 1}\end{pmatrix}
\]
satisfy the conditions 
\refeqs(C-1,C-2) with for $m=k$ as $M(n(Nd+1))$-valued functions. 
Hence, the induction hypothese implies 
\[
\sup_{t\in \R}\|{\textbf u}(t)\|_{\Hs^k} 
\leq C \|{\textbf u}(0)\|_{\Hs^{k}}, 
\]
which implies \refeq(0-3aa) and \refth(1) 
is proved for all integer s $k\geq 0$.  
\qed 

\section{Proof of \refthb(2)} \lbsec(proof-2)
We apply Kato's theory of smooth perturbation (\cite{Kato-e}). 
Let $\tilde{\Hs}$ and $\tilde{\Hs}'$ be separable Hilbert spaces 
(possibly $\tilde{\Hs}=\tilde{\Hs}'$), $T$ a selfadjoint operator on 
$\tilde{\Hs}$ and $A$ a closed operator from $\tilde{\Hs}$ to $\tilde{\Hs}'$. 
$A$ is said to be $T$-smooth if $A$ is $T$-bounded and 
for any $u \in \tilde{\Hs}$ 
\bqn \lbeq(Kato-s)
\sup_{\ep>0} \int_{\R} \|A (T-\lam \pm i\ep)^{-1} u\|_{\tilde{\Hs}'}^2 d\lam 
= \frac1{2\pi}\int_0^\infty \|A e^{\pm it T} u\|_{\tilde{\Hs}'}^2 dt 
\leq C \|u\|^2.  .   
\eqn 

Recall $\Kg_0$ is the selfadjoint operator on $\Ks= L^2(\R) \otimes \Hs$ 
defined by \refeq(U0). We have the obvious identity:
\bqn \lbeq(Phillips-0)
(\Kg_0- \z)^{-1} = \pm i \int_0^\infty e^{\mp (\Kg_0-\z)\s }d\s, \quad 
\z \in \C^{\pm} .
\eqn
Let $A_{jk}(t)$ and $B_{jk}(t)$ be the multiplication operators on $\Hs$ with 
$A_{jk}(t,x^{jk})$ and $B_{jk}(t,x^{jk})$ respectively and 
$\Ag_{jk}$ and $\Bg_{jk}$ be multiplication operators on $\Ks$ defined by 
\[
\Ag_{jk} u(t,x) = A_{jk}(t)u(t,x), \quad 
\Bg_{jk} u(t,x) = B_{jk}(t)u(t,x)\,.
\]

\bglm  \lblm(1) 
Suppose that the conditions of \refth(2) are satisfied. 
Then $\Ag_{jk}$ and $\Bg_{jk}$ are $\Kg_0$-smooth. 
\edlm 
\bgpf We prove the lemma for $\Ag_{jk}$. The proof for $\Bg_{jk}$ is similar. 
We write $\Sg$ for $\Sg(\R^{1}\times X\colon \C^n)$ for simplicity. 
For $u \in \Sg$, it is easy to see that  
$(e^{-i\s\Kg_0}u)(t,\cdot)= e^{-i\s{H}_0} u(t-\s, \cdot)\in \Sg$  
and $\Ag_{jk}e^{-i\s{\Kg_0}}u \in \Ks$. Then, changing the variable 
$t $ to $t+\s$ and changing the order of integrations, we have  
\bqn \lbeq(CV)
\int_{\R} \|\Ag_{jk}e^{-i\s{\Kg_0}}u\|_{\Ks}^2 d\s
= \int_{\R} 
\left( \int_{\R}\|A_{jk}(t+\s)e^{-i\s{H_0}}u(t)\|_{\Hs}^2 d\s \right)dt.
\eqn 
Since $e^{-i{\s}H_0 }= e^{i{\s}\lap_{x^{jk}}} \otimes e^{i{\s}\lap_{x_{jk}}}$ 
and since $e^{i{\s}\lap_{x_{jk}}}$ is unitary and commutes with 
$A_{jk}(t,x^{jk})$, H\"older's inequality 
with respect to the variable $x^{jk}$ implies 
\begin{align}
& \|A_{jk}(t+\s)e^{-i\s{H_0}}u(t)\|_{\Hs}^2 
= \int_{X_{jk}} 
\|A_{jk}(t+\s)e^{i\s\lap_{x^{jk}}}u(t, x^{jk}, x_{jk})\|_
{L^2(X_{jk}^r)}^2 dx_{jk} \notag \\
& \leq 
\sup_{t\in \R}\|A_{jk}(t)\|_{L^d}^2 \int_{X_{jk}^c} \|
e^{i\s\lap_{x_{jk}}}u(t,x^{jk},{x_{jk}})\|_{L^{\frac{2d}{d-2}}(X_{jk}^r)}^2 
dx_{jk}  \lbeq(Holder)
\end{align}
We denote $C_{jk}=\sup_{t\in \R}\|A_{jk}(t)\|_{L^d}^2$ and integrate 
both sides of \refeq(Holder) by $d\s$. Changing the order of integrations 
and applying the end point Strichartz inequality (\cite{KT}):
\[
\int_{\R}
\|e^{i\s\lap_{x^{jk}}}
u(t,x^{jk},x_{jk})\|_{L^{\frac{2d}{d-2}}(X_{jk}^r)}^2 d\s  
\leq C_{KT} \|u(t,x^{jk},x_{jk})\|_{L^{2}(X_{jk}^r)}^2, 
\]
we have that there exists a constant $C$ independent of $u$ and $t$ such that  
\bqn \lbeq(ext)
\int_0^\infty \|A_{jk}(t+\s)e^{-i\s{H_0}}u(t)\|_{\Hs}^2 d\s 
\leq C\|u(t)\|_{\Hs}^2.  
\quad u \in \Sg. 
\eqn 
On substituting the result into \refeq(CV) we obtain 
\bqn \lbeq(2-1) 
\int_{\R} \|\Ag_{jk}e^{-i\s{\Kg_0}}u\|_{\Ks}^2 d\s \leq C \|u\|_{\Ks}^2, 
\ \ u \in \Sg. 
\eqn  
It follows by Schwarz' equality that for $\ep>0$  
\bqn \lbeq(2-3)
\int_0^\infty e^{-\s \ep} \|\Ag_{jk}e^{-i\s\Kg_0}u \|_{\Ks} d\s 
\leq C (2\ep)^{-\frac12}  \|u\|_{\Ks}.
\eqn 
Then, by applying Minkowski's inequality, we obtain from \refeq(2-3) that 
\bqn \lbeq(2-0)
\|\Ag_{jk}(\Kg_0-\lam -i\ep)^{-1} u \|_{\Ks}
= \left\| \int_0^\infty e^{i\s\lam } e^{-\s \ep}\Ag_{jk}e^{-i\s\Kg_0}u d\s
\right\|_{\Ks}
\leq C \ep^{-\frac12}  \|u\|_{\Ks}.
\eqn 
Since $\Ag_{jk}$ is closed and since $\Sg$ is dense in $\Ks$, both of 
\refeq(2-1) and \refeq(2-0) extend to all $u \in \Ks$ and 
$\Ag_{jk}$ is $\Kg_0$-smooth. 
\edpf 

Let $\Ks_N = \Ks \otimes \C^{N(N-1)/2}= \Ks^{\oplus(N(N-1)/2)}$ 
and $\Ag$ and $\Bg$ be the multiplication operators from 
$\Ks$ to  $\Ks_N $ defined by 
\[
\Ag u  = \oplus_{j<k} \Ag_{jk}u, \quad 
\Bg u = \oplus_{j<k} \Bg_{jk} u
\]
respectively so that the multiplication operator with 
$V(t,x)$ on $\Ks$ is equal to $\Bg^\ast \Ag$. Denote 
$\Rg_0(\z)= (\Kg_0 - \z)^{-1}$. 
 
\bglm \lblm(3) 
\ben 
\item[\textrm{(1)}] The operators $\Ag$ and $\Bg$ are $\Kg_0$-smooth. 
\item[\textrm{(2)}] There exists a constant $C$ independent of 
$\z \in \C^{\pm}=\{\pm \Im z> 0\}$ such that  
\bqn \lbeq(smooth-resol)
\|\Ag \Rg_0(\z) \Bg^\ast u\|_{\Ks_N }\leq C \|u\|_{\Ks_N} , \quad 
u \in D(\Bg^\ast).
\eqn 
\een
\edlm 
\bgpf Statement (1) is obvious by \reflm(1) and we prove (2). 
The following is an adaptation of 
Iorio-O'Carrol's argument (\cite{IO}) to the 
Howland scheme. It suffices to show that 
\bqn \lbeq(20)
\sup_{\pm \Im \z >0}
\|\Ag_{jk} \Rg_0(\z) \Bg_{lm}^\ast u\|_{\Ks}\leq C \|u\|_{\Ks}, \quad 
j<k, l<m.
\eqn  
We prove the $+$ case. The proof for the other case is similar. 
For $u,v \in \Sg$, we have  as in \refeq(2-0) that 
\bqn \lbeq(3-2)
|(\Rg_0(\z) \Bg_{lm}^\ast u, \Ag_{jk}^\ast v)_{\Ks}|
\leq 
\int_0^\infty |(e^{-i\s\Kg_0}\Bg_{lm}^\ast u, \Ag_{jk}^\ast v)_{\Ks}| d\s.
\eqn

\noindent 
1) Let first $\{j,k\} \cap \{l,m\}=\emptyset$. Define $X_{jklm}$ by 
$X= X_{jk}^r \oplus X_{lm}^r \oplus X_{jklm}$ and let 
$x=x^{jk}\oplus x^{lm} \oplus x'$ be the corresponding orthogonal 
decomposition. Then 
\[
H_0= - \lap_{jk} - \lap_{lm} - \lap', \ \ 
\lap_{jk} \colon =\lap_{x^{jk}}, \ \lap_{lm}\colon =\lap_{x^{lm}}, \ 
\lap' \colon =\lap_{x'}.
\]
Since 
$e^{i\s(\lap_{jk}+\lap')}$ and 
$e^{i\s\lap_{lm}}$ commute with the multiplications 
$B_{lm}^\ast(t)$ and  
$A_{jk}^\ast(t)$ respectively, 
$(e^{-i\s\Kg_0}\Bg_{lm}^\ast u, \Ag_{jk}^\ast v)_{\Ks}$ 
is equal to 
\begin{align}
& \int_{\R} (e^{i\s\left(\lap_{jk}+ \lap_{lm}+ \lap'\right)}
B_{lm}^\ast(t-\s)u(t-\s), A_{jk}^\ast(t)v(t))_{\Hs} dt \notag \\
& = \int_{\R} (
A_{jk}(t) e^{i\s\left(\lap_{jk}+ \lap'\right)}u(t-\s), 
B_{lm}(t-\s)e^{-i\s\lap_{lm}} v(t))_{\Hs} dt. \lbeq(21)
\end{align}
Then, since $e^{i\s\lap'}$ commutes with $A_{jk}(t)$ and is unitary, 
Schwarz's inequality implies 
\begin{align}
|(e^{-i\s\Kg_0}\Bg_{lm}^\ast u, \Ag_{jk}^\ast v)_{\Ks}| & \leq  
\left(\int_{\R} 
\|A_{jk}(t) e^{i\s\lap_{jk}}u(t-\s)\|_{\Hs}^2 dt \right)^\frac12 
\notag \\
& \times \left(\int_{\R} 
\|B_{lm}(t-\s) e^{i\s\lap_{lm}} v(t)\|_{\Hs}^2 dt\right)^\frac12.
\lbeq(22)
\end{align} 
Change variable $t$ by $t+\s$ in the first factor, 
integrate both sides of \refeq(22) by $d\s$, 
apply Schwarz's inequality once more 
and change the order of integrartions. We obtain   
\begin{align} 
& \int_{0}^\infty |(e^{-i\s\Kg_0}\Bg_{lm}^\ast u, \Ag_{jk}^\ast v)_{\Ks}|d\s
\leq 
\left(\int_{\R} \left(\int_{0}^\infty 
\|A_{jk}(t+\s) e^{i\s\lap_{jk}}u(t)\|_{\Hs}^2 d\s\right)dt 
\right)^{\frac12} \notag \\
& \hspace{2cm} 
\times 
\left(\int_{\R} \left(\int_{0}^\infty  
\|B_{lm}(t-\s) e^{i\s\lap_{lm}}v(t)\|_{\Hs}^2 d\s\right)dt 
\right)^{\frac12} \lbeq(23)
\end{align}  
Then, repeating the argument used for \refeq(2-3), 
we obtain by using H\"older's, the end point Strichartz inequality 
and the density of $\Sg$ in $\Ks$ that 
\bqn \lbeq(24)
\refeq(23) \leq C \sup_{t\in \R}\|A_{jk}(t,\cdot)\|_{L^d}
\sup_{t \in \R}\|B_{lm}(t,\cdot)\|_{L^d}\|u\|_{\Ks} \|v\|_{\Ks}\,, \ 
u, v \in \Ks.
\eqn 
It follows from \refeq(3-2) and \refeq(24) that  there exists a 
constant $C>0$ independent of $u, v \in \Ks$ and $\z\in \C^{+}$ such that 
\[
|(\Rg_0(\z) \Bg_{lm}^\ast u, \Ag_{jk}^\ast v)_{\Ks}| \leq 
C\|u\|_{\Ks} \|v\|_{\Ks}
\]
and we obtain the desired \refeq(20) for the case 
$\{j,k\} \cap \{l,m\}=\emptyset$. 

2) Let next $\{j,k\}=\{l,m\}$. Let $u \in \Sg$ and $\lap_{jk}=\lap_{x^{jk}}$. 
We sandwitch \refeq(Phillips-0) by $\Ag_{jk}$ and $\Bg_{jk}$ and apply 
Minowski's inequlaity. We obtain for any $\z\in \C^{+}$ that 
\begin{align} 
& \|\Ag_{jk} \Rg_0(\z) \Bg_{jk}^\ast u \|_{\Ks}
\leq 
\int_0^\infty \|\Ag_{jk}e^{-i\s\Kg_0}\Bg_{jk}^\ast u \|_{\Ks}d\s \notag \\
& \leq 
\int_0^\infty \left(\int_{\R}
\|A_{jk}(t)e^{i\s\lap_{jk}} B_{jk}^\ast(t-\s) u(t-\s) \|_
{L^2(X_{jk}^r)\otimes L^2(X_{jk})}^2 dt \right)^{\frac12} d\s. \lbeq(3-3)
\end{align}
The well-known $L^p$-$L^q$ estimates for $e^{i\s\lap_{jk}}$ on $L^2(X_{jk}^r)$ 
and H\"older's inequality imply that the integrand of the inner  
integral of \refeq(3-3) is bounded for any $\r\geq 2$ by 
\bqn \lbeq(25)
C_r 
(\sup_{t\in \R}\|A_{jk}(t,\cdot)\|_{L^\r(X_{jk}^r)} 
\sup_{t\in \R}\|B_{jk}(t,\cdot)\|_{L^\r(X_{jk}^r)})
|\s|^{-\frac{d}{\r}} \|u(t-\s)\|_{\Hs}^2.
\eqn 
Since $A_{jk}(t,\cdot)$ and $B_{jk}(t,\cdot)$ are 
$(L^p\cap L^q)(X_{jk}^r)$-valued 
bounded function of $t\in \R$ for 
$p,q$ such that $2\leq p <d<q\leq \infty$, we have by applying \refeq(25) 
for $\r=p$ and $\r=q$ that 
\[
\refeq(3-3) \leq C \int_0^\infty \min_{r=p,q} |\s|^{-\frac{d}{r}} 
\left(\int_{\R}\|u(t-\s)\|^2_{\Hs}dt\right)^{\frac12}d\s \leq 
C \|u\|_{\Ks}.
\]
Hence \refeq(20) holds also for $\{j,k\}=\{l,m\}$. \\[2pt]

\noindent 
3) Finally let 
$|\{j,k\} \cap \{l,m\}|=1$. We may assume $N=3$ and 
$\{j,k\}=\{1,2\}$ and $\{l,m\}=\{2,3\}$ without losing generality. 
We use Jacobi-coordinates for three particle:  
\[
r=x_2-x_1, \ \  y = x_3 - \frac{m_1 x_1 + m_2 x_2}{m_1+ m_2}, \ \ 
x_{cm}= \frac{m_1 x_1+ m_2 x_2 + m_3 x_3m }{m_1+ m_2+ m_3}.
\]
Let $\m$, $\n$ be reduced masses and $m$ the total mass:
\[
\m=m_1m_2/(m_1+m_2), \ \ \n=(m_1+m_2)m_3/(m_1+m_2+m_3,\ \ 
m= m_1+ m_2+ m_3.
\]
Then, $x_3-x_2= y-ar$, $a=m_1 r/(m_1+m_2)$, $dx = Cdr dy dx_{cm}$ and 
\[
\lap_{x}= \frac{1}{2\m}\lap_{r}+\frac1{2\n}\lap_{y}+ \frac1{2m}\lap_{x_{cm}}.
\]
It follows, by writing $u(x)= u(r,y,x_{cm})$, that  
\begin{multline} \lbeq(jacobi)
(\Ag_{12}e^{-i\s\Kg_0}\Bg_{23}^\ast u)(t,r,y,x_{cm}) \\
= A_{12}(t,r)
e^{i\s\left(\frac{\lap_r}{2\m}+ \frac{\lap_y}{2\n}+ 
\frac{\lap_{x_{cm}}}{2m}\right)} 
B_{23}^\ast (t-\s, ar+y)u(t-\s,r,y,x_{cm})\, .
\end{multline} 
Since $e^{i\s(\lap_{y}/2\n+ \lap_{x_{cm}}/2m)}$ commutes with $A_{12}(t,r)$ 
and is unitary on $\Hs$, \refeq(jacobi) implies that 
$\|\Ag_{12}e^{-i\s\Kg_0}\Bg_{23}^\ast u\|_{\Ks}^2$ is equal to
\[
\int_{X_{c} \times \R}
\|A_{12}(t,r)e^{i\s\lap_{r}/2\m} 
B_{23}^\ast (t-\s, ar+y)u(t-\s,r,y,x_{cm})\|_{L^2(dr)}^2 dy dx_{cm} dt. 
\]
We fix $(t,y,x_{cm})$, apply $L^p$-$L^q$ 
estimates for $e^{i\s\lap_{r}/2\m}$ and H\"order's inequality 
on $L^2(dr)$ as in the case 2) and estimate the integrand by 
\[
C 
\min_{\r\in \{p,q\}} \{\|A_{12}(t)\|_{\r} \|B_{23}^\ast (t-\s)\|_{\r}
|\s|^{-d/\r}\}
\|u(t-s,r,y,x_{cm})\|_{L^2(dr)}.
\]
This yields that $\|\Ag_{jk}e^{-i\s\Kg_0}\Bg_{kl}^\ast u\|_{\Ks}\leq 
C \min_{r\in \{p,q\}} |\s|^{-d/\r} \|u\|_{\Ks}$ and we obtain as previously 
the desired estimae \refeq(20) for this case. 
This completes the proof of \reflm(3). 
\edpf 

It follows from \reflm(3) that Theorems 1.5 ans 3.9 of \cite{Kato-e} may be 
applied to the triplet $(\Ag, \Bg, \Kg_0)$ on $\Ks$ and we obtain 
the following theorem. 

\bgth[Kato(\cite{Kato-e})] \lbth(kato)
There exists a constant $\ka_0>0$ such that for $|\ka|<\ka_0$ 
there is a unique closed operator $\tKg(\ka)$ on $\Ks$ which satisfies 
following such that {\textrm (a)} to {\textrm (d)}. $\tKg(\ka)$ 
is uniquely determined by {\textrm (a)} and {\textrm (b)}. 
\ben 
\item[{\textrm (a)}] $\Ag$ is $\tKg(\ka)$-smooth and $\Bg $ is 
$\tKg^\ast(\ka)$-smooth.
\item[{\textrm (b)}] The resolvent 
$\tilde{\Rg}(\z, \ka) = (\tKg(\ka) -\z)^{-1}$ exists for all 
$\z\in \C \setminus \R$ and itsatisfies 
\bqn \lbeq(RR0)
\tilde{\Rg}(\z, \ka)- {\Rg}_0(\z)
=- \ka [\Rg_0(\z)\Bg^\ast]\Ag \tilde{\Rg}(\z, \ka)
= -\ka [\Rg(\z,\ka){\Bg}^\ast]{\Ag}\Rg_0(\z)
\eqn 
where $[\Rg_0(\z){\Bg}^\ast]$ is the closure of $\Rg_0(\z){\Bg}^\ast$ 
and similarly for $[\Rg(\z,\ka){\Bg}^\ast]$. 
\een
Moreover, following statements are satisfied: 
\ben 
\item[{\textrm (c)}] The equations 
\[
(\Wg_\pm(k) u,v)= (u,v)\mp \lim_{\ep \to 0}
\frac{1}{2{\pi}i}
\int_{\R} (\Ag \Rg_0(\lam\pm i\ep)u, 
\Bg \tilde{\Rg}(\lam\mp i\ep,\ka)^\ast v)d\lam 
\]
define bounded operators in $\Ks$ which have bounded inverses 
and $\tKg(\ka)$ and $\Kg_0$ are similar to each other via 
$\Wg_{\pm}(\ka)$ and $\Wg_{\pm}(\ka)^{-1}$:  
\bqn \lbeq(similar)
\tKg(\ka) = \Wg_\pm (\ka)\Kg_0 \Wg_\pm(\ka)^{-1}. 
\eqn 
\item[{\textrm (d)}] $\tKg(\ka)$ generates 
a strongly continuous group of uniformly bounded operators 
$ e^{-i\s \tKg(\ka)}$ in $\Ks$: $\| e^{-i\s \tKg(\ka)}\|_{\Bb(\Ks)}\leq M$ 
and 
\bqn \lbeq(6-2)
\Wg_\pm (\ka)= \lim_{\s \to \pm \infty} e^{i\s\tKg(\ka)} e^{-i\s\Kg_0}, \quad 
\Wg_\pm(\ka)^{-1}
= \lim_{\s \to \pm \infty} e^{i\s\Kg_0} e^{-i\s\tKg(\ka)}\,.
\eqn 
\een 
\edth 

In view of \refth(kato), the next lemma completes the proof of \refth(2). 
Recall that $\Kg(\ka)$ is the generator of $\Ug(\s,\ka)$,see  \refeq(generator).
\bglm We have the equality $\Kg(\ka)= \tKg(\ka)$, $\tKg(\ka)$ being as in 
\refth(kato). 
\edlm 
\bgpf Let $M$ be as in \refeq(bdd) and $\tM=M+1$ and 
$\Rg(z,\ka)= (\Kg(\ka)-\z)^{-1}$ for $|\Im \z |>\tM$ (see \refeq(Phillips)). 
It suffices to show     
\bqn  
{\Rg}(\z, \ka)- \Rg_0(\z)
= -\ka [\Rg_0(\z)\Bg^\ast] {\Ag}{\Rg}(\z, \ka), \ \ |\Im \z |>\tM. \lbeq(resol)
\eqn 
Indeed \refeq(resol) implies $\Ag{\Rg}(\z, \ka)- \Ag\Rg_0(\z)
= -\ka \Ag [\Rg_0(\z)\Bg^\ast] {\Ag}{\Rg}(\z, \ka)$ by virtue of 
\reflms(1,3) and 
\bqn \lbeq(RR0-a)
{\Rg}(\z, \ka)= {\Rg}_0(\z, \ka)
- \ka [\Rg_0(\z)\Bg^\ast]
(1+ \ka [\Ag {\Rg}_0(\z){\Bg}^\ast])^{-1}{\Ag}\Rg_0(\z) \,.
\eqn 
Hence ${\Rg}(\z, \ka)=\tilde{\Rg}(\z, \ka)$ for $|\Im \z |>\tM$ and 
$\Kg= \tKg$. 

We now prove \refeq(resol). We prove it for $\Im \z > \tM$. The proof for 
$\Im \z<-\tM$ is similar. In what follows we omit the coupling constant $\ka$. 
We first show that 
\bqn \lbeq(inte)
\int_{0}^\infty  e^{-\tM \s} \|\Ag \Ug(\s)u\|_{\Ks_N} d\s \leq C \|u\|_{\Ks}.
\eqn 
Let $\Xg_{jk}=L^{\frac{2d}{d-2}}(X_{jk}^r,\Hs_{jk})$ 
and $C_d=\sum_{j<k} \sup_{t\in \R}\|A_{jk}(t,\cdot)\|_{L^d(X_{jk}^r)}$. Then, 
by H\"older's inequality 
\begin{align*}
& \int_0^\infty \|\Ag \Ug(\s)u\|_{\Ks_N}^2 e^{-2\tM\s}d\s \\
& = \sum_{j<k}\int_{\R}
\left(\int_{0}^\infty 
e^{-2\tM\s}\|A_{jk}(t+\s)U(t+\s,t)u(t)\|_{\Hs}^2 d\s\right) dt \\
& \leq C_d^2  
\int_{\R}
\left(\int_{0}^{\infty} e^{-2\tM\s}\|U(t+\s,t)u(t)\|_{\Xg_{jk}}^2 d\s\right) 
dt\,. 
\end{align*} 
The inner $d\s$-integral is equal to 
\begin{align}
& \sum_{n=0}^\infty 
\int_{0}^1 e^{-2\tM(\s+n)}\|U(t+\s+n,t)u(t)\|_{\Xg_{jk}}^2 d\s \notag \\
& = \sum_{n=0}^\infty 
\int_{0}^1 e^{-2\tM(\s+n)}\|U(t+\s+n,t+n)U(t+n, t)u(t)\|_{\Xg_{jk}}^2 d\s.
\lbeq(ds)
\end{align}
The time local Strichartz estimate \refeq(local-Strichartz) with 
$p=d/2$  implies that 
\[
\sup_{t\in \R,n=0,1,\dots}
\int_0^{1}\|U(t+\s+n,t+n,\ka)\ph\|_{\Xg_{jk}}^2 d\s 
\leq C \|\ph\|_{\Hs}^2
\]
and, the right of \refeq(ds) is bounded by 
\[
C \sum_{n=0}^\infty e^{-2\tM{n}}\|U(t+n, t)u(t)\|_{\Hs}^2 .
\]
with a constant independent of $t$. Combining these estimates yields  
\begin{align*}
& \int_0^\infty \|\Ag \Ug(\s)u\|_{\Ks_N}^2 e^{-2\tM\s}d\s 
\leq C C_d^2 
\sum_{n=0}^\infty \int_{\R}e^{-2\tM{n}}\|U(t+n, t)u(t)\|_{\Hs}^2 dt \\ 
& \hspace{2cm} 
= CC_d^2 \sum_{n=0}^\infty e^{-2\tM{n}}\|\Ug(n)u\|_{\Ks}^2
\leq C \sum_{n=0}^\infty e^{-2n} \|u\|_{\Ks}^2 
= C \|u\|_{\Ks}^2.
\end{align*}
We likewise have the corresponding estimate for the integral on 
$\s \in (-\infty,0)$. Thus, the Schwarz' inequality implies \refeq(inte) 
and 
\bqn \lbeq(Resol)
\int_{0}^\infty e^{\mp {i\z\s}} \Ag \,\Ug(\pm \s)u d\s = \Ag\Rg(\z) u, 
\quad u \in \Ks
\eqn 
is uniformly bounded for $|\Im \z|>\tM$:  
\bqn 
\|\Ag \Rg(\z)u\|_{\Ks} \leq C \|u\|_{\Ks}, \quad |\Im \z|>\tM. 
\eqn 
As a function of $t\in \R$, 
$U(t,s)\ph\in C(\R, \Hs) \cap \Xg_{\loc}^{d/2}$ and $\Ag$ and $\Bg$ map  
$C(\R,\Hs) \cap \Xg_{\loc}^{d/2}$ into $L^2_{\textrm{loc}}(\R, \Hs)$ 
and the contrction of the solution (\cite{Ya-2016}) implies the Duhamel formula 
\bqn \lbeq(-1)
U(t,s)\ph= 
e^{-i(t-s)H_0}\p - i\int_s^t e^{-i(t-r)H_0}B^\ast(r) A(r)U(r,s)\ph dr .
\eqn 
Let $u,v \in \Sg$. 
Substituting $t-\s$ for $s$ and $u(t-\s)$ for $\ph$ in \refeq(-1) yields 
\begin{multline*} 
U(t,t-\s)u(t-\s)= 
e^{-i{\s}H_0}u(t-\s) \\ 
- i\int_0^\s e^{-i(\s-\t)H_0}B^\ast(t-\s+\t)
A(t-\s+\t)U(t-\s+\t, t-\s)u(t-\s)d\t,
\end{multline*} 
which implies that 
\[
(\Ug(\s)u)(t) = (\Ug_0(\s)u)(t)- 
i \int_{0}^\s (\Ug_0(\s-\t)\Bg^\ast \Ag \Ug(\t)u)(t) d\t, \ \ t \in \R.
\]
Taking the inner products with $v \in \Sg$ in both sides produces 
\[
(\Ug(\s)u, v)_{\Ks}= (\Ug_0(\s)u, v)_{\Ks} 
- i \int_{0}^\s (\Ag \Ug(\t)u, \Bg \Ug_0(\t-\s)v)_{\Ks} d\t\,.
\]
Multiply both sides by $i e^{i\s\z}$ with $\Im \z>\tM$ 
and integrate by $d\s$ over $[0,\infty)$. We obtain by changing 
the order of integrations that 
\[
(\Rg(\z)u,v)_{\Ks}= (\Rg_0(\z)u,v)_{\Ks} - (\Ag \Rg(\z)u, 
\Bg \Rg_0(\overline{\z})v)_{\Ks}.
\]
Restoring the coupling constant, we obtain \refeq(resol): 
for $\Im z>\tM$. This completes the proof. 
\edpf


\end{document}